\documentclass[11pt,a4paper,reqno]{amsart}
\usepackage{graphics}
\usepackage{amsmath}

\subjclass[2000]{Primary: 37E05, 37B10 
Secondary: 37A05, 37B40}
\keywords{interval translation map, interval exchange transformation,
substitution shift, unique ergodicity, Hausdorff dimension}

\newtheorem{theorem}{Theorem}
\newtheorem{proposition}[theorem]{Proposition}
\newtheorem{corollary}[theorem]{Corollary}
\newtheorem{lemma}[theorem]{Lemma}
\renewenvironment{proof}{\noindent {\bf Proof.}}{ \hfill\qed\\ }
\newenvironment{proofof}[1]{\noindent {\bf Proof of #1.}}{ \hfill\qed\\ }

\def\hT{\hat{T}}
\def\K{\mathcal{K}}
\def\U{\mathcal{U}}
\def\V{\mathcal{V}}
\def\e{\varepsilon}
\def\Ae{A_{\e}}\def\Fe{F_{\e}}
\def\Ne{N_{\e}}
\def\O{\Omega}
\def\Ob{\overline{\O}}
\def\a{\alpha}
\def\b{\beta}
\def\ab{(\a,\b)}
\def\Tab{T_{\a,\b}}
\def\Sab{\Sigma_{\a,\b}}
\def\freq{\mbox{freq}}
\def\CH{\mbox{CH}}
\def\N{\mathbb{N}}

\def\eps{\varepsilon}
\def\RR{{\mathbb R}}

\def\NN{{\mathbb N}}
\def\OO{{\mathcal O}}
\def\1{{{\mathit 1} \!\!\>\!\! I} }
\renewcommand{\phi}{\varphi}
\newcommand{\inte}{\mathop{int}}

\def\diam{{\hbox{{\rm diam}}\,}}
\newcommand{\eqdef}{\stackrel{\scriptscriptstyle\rm def}{=}}

\def\orb{\hbox{{\rm orb}}}

\def\ubd{\overline{\mbox{dim}_{\rm B}}}
\def\lbd{\underline{\mbox{dim}_{\rm B}}}
\def\hd{\mbox{dim}_{\rm H}}
\def\rl{{\underline r}}
\def\ru{{\overline r}}

\begin{document}
\bibliographystyle{plain}
\title{The Gauss map on a class of interval translation mappings}
\author{H.~Bruin \and S.~Troubetzkoy}
\thanks{We gratefully acknowledge the support of 
Universit\'e de Toulon et du Var,
the Royal Netherlands Academy of Arts and Sciences (KNAW)
and Van Gogh NWO/EGIDE
}
\address{Department of Mathematics,
University of Groningen,
P.O. Box 800, 9700 AV Groningen, the Netherlands}
\email{bruin@math.rug.nl}
\urladdr{http://www.math.rug.nl/{\lower.7ex\hbox{\~{}}}bruin/}

\address{Centre de Physique Th\'eorique, Institut de Math\'ematiques de
Luminy, F\'ed\'eration de Recherche des Unit\'es de Math\'ematiques de 
Marseille and Universit\'e de la M\'editerran\'ee,  Luminy, Case 907,
F-13288 Marseille Cedex 9, France}   
\email{troubetz@iml.univ-mrs.fr}
\urladdr{http://iml.univ-mrs.fr/{\lower.7ex\hbox{\~{}}}troubetz/}

\begin{abstract}
We study the dynamics of a class of interval translation map on three
intervals. We show that in this class the typical ITM is of finite type
(reduce to an interval exchange transformation)
and that the complement  
contains a Cantor set.  We relate
our maps to substitution subshifts.
Results on Hausdorff dimension of the attractor and on unique 
ergodicity are obtained.
\end{abstract}
\maketitle

\section{Introduction}
In the last years there has been an increased interest in the study of piecewise
isometries, see e.g.~ \cite{AKT,ACP,G,L,LV,M,V1,V2} and the references 
therein.  Classically, invertible
piecewise isometries have been studied in one dimension;
these are the so-called interval exchange mappings (IETs for short), which
can appear for example as first return maps of geodesic and billiard flows.  
A natural generalization of IETs to the noninvertible case,
interval translation mappings (ITMs for short), has been recently introduced
by Boshernitzan and Kornfeld \cite{BK}. 
ITMs have been studied for their topological dynamics in \cite{ST}, 
and for their invariant measures and complexity in \cite{BH}. 
Many 0-entropy maps of the interval are (semi)conjugate to ITMs
\cite{BH}.
In \cite{BK} several very interesting questions were asked,
which we answer in this article for a special class of interval
translation mappings.  

We define a class of ITMs which can be viewed as
translations of two intervals on the circle.
The example given in \cite{BK} is a special member of our class. 
This example was the first one for which a Cantor attractor for an ITMs was
observed.
We define an inducing procedure similar to Rauzy induction for interval
exchange mappings \cite{R,V1}. It defines a map $G$ in parameter space which 
plays the same
role as the well-known Gauss map for circle rotations; 
we call $G$ the {\em Gauss map} for our class of ITMs.     
Whether the induction procedure can be extended to all
ITMs is a very interesting open problem.

Using the Gauss map, we prove that almost
every map in our class is of finite type (i.e., its attractor is a
union interval rather than a Cantor set), giving a partial answer to a
question posed in \cite{BK}. If the ITM has a Cantor attractor, we specify
an isomorphism to a shift space generated by a chain of substitutions.
We give an upper bound
on the dimension of the attractor. Although the attractor is ``dynamically defined'',
it is interesting to note that its
upper box dimension need not be equal to its Hausdorff dimension.
Hausdorff measure, whether finite or not, is always $T$-invariant. 
Finally we give sufficient conditions for mappings in
our class to be uniquely ergodic, and a sufficient condition preventing
unique ergodicity. 

In a companion paper, J.~Cassaigne shows that the (subword) complexity $p(n)$ of the
subshift describing our system is linear, in fact $p(n) \leq 3n$ \cite{Cas2003}.
This yields a partial answer to another question posed in \cite{BK}. 
Cassaigne's original technique \cite{Cas94b} is designed for substitution subshifts.
In \cite{Cas2003}, he extends the method for the chain of substitutions developed in our paper.

\section{Statement of results}
Let $0 = \b_0 < \b_1 < \dots < \b_r = 1$, $I= [0,1)$ and, for $i=0, \dots, r$,
$B_i \eqdef [\b_{i-1}, \b_i)$.
An interval translation mapping is an interval map
$T:I \to I$ given by
\[
T(x) \eqdef x+\gamma_i \mbox{ if } x \in B_i,
\]
where $\gamma_i$ are fixed numbers such that $T$ maps $I$
into itself.
We also define the image of $1$ by $T(1) \eqdef \lim_{x \to 1^-} T(x)$.

Define $\Omega_0 = I$ and $\Omega_n = T(\Omega_{n-1})$.
The set $\Omega \eqdef {\cap_n \Omega_n}$ 
is called the {\em attractor} of the ITM.
We say that $T$ is of {\em finite type} if 
$\Omega_n = \Omega_{n+1}$ for some $n$. In this case, the
attractor $\Omega = \Omega_n$ is a finite union of intervals, and
$T|_{\Omega}$ is an interval exchange transformation (IET).
If $\Omega_{n+1}$ is strictly smaller than $\Omega_n$ for each $n$,
then $\overline{\Omega}$ is a Cantor set or the union of a Cantor set
and a finite collection of intervals.
The latter only happens when $T$ is reducible; orbits do not visit both the
Cantor set and the intervals. We will ignore this reducible case, and 
concentrate on the case that  
$\overline{\Omega}$ is a Cantor set; $T$ is said to be
of {\em type} $\infty$ in this case. For convenience, we will often consider
$\overline{\Omega}$ instead of $\Omega$ to be the attractor, and assume that $T|_{\overline{\Omega}
\setminus \Omega}$ is (re)defined by continuity from the left.

Our first result is a general structure theorem about the topological
dynamics of ITMs. In fact, this theorem is a consequence of results from 
\cite{HR}, but since that proof requires extensive machinery and is 
spread out over several papers, we prefer to give a direct proof for our case,
which is more in the spirit of the well-known results for
IETs, see \cite{KH}.

\begin{theorem}\label{thmTrans}
If $T|_{\overline{\Omega}}$ is transitive, then this restriction is minimal.
\end{theorem}

It follows from \cite{ST} that $\Ob$ is a Cantor set when $T$ is of infinite 
type. We have 

\begin{theorem}\label{thmHMeas}
Let $d$ be the Hausdorff dimension of ${\overline{\Omega}}$. 
The Hausdorff
$d$ dimensional measure  $H_d$ on ${\overline{\Omega}}$ is $T$-invariant.
\end{theorem}

Remark: this theorem is a new result for ITMs of infinite type.  
In general we do not know if $0 < H_d(\Ob) < \infty.$ In fact, if $d=0$ then
$H_d(\Ob)$ is infinite. Examples of ITMs for which $d=0$ are produced
in Theorem \ref{thmHD}.

The main results of this article are on 
interval translation maps of a special form.
Consider $U \eqdef \{ \ab : 0\le \b \le \a \le 1\}$,
$L \eqdef \{ \ab : 0 \le \a \le \b + 1  \le 1\}$
and $R \eqdef U \cup L$. 
For $\ab$ in the interior $U^{\circ}$ of $U$ consider the ITM
$T = T_{\a,\b}:[0,1) \to [0,1)$ defined by (see Figure 1)
\[
T(x) \eqdef \left \{
\begin{array}{ll} 
x + \a  & \text{ for } x \in [0,1-\a)\cr
x + \b  & \text{ for } x \in [1-\a,1-\b)\cr
x + \b-1 & \text{ for } x \in [1-\b,1). 
\end{array}\right.
\]
\begin{figure}[ht]
\unitlength=7mm
\begin{picture}(8,8)
\thinlines
\put(1,1){\line(1,0){6}}
\put(1,7){\line(1,0){6}}
\put(7,1){\line(0,1){6}}
\put(1,1){\line(0,1){6}}

\put(7.2,1.1){\vector(0,1){0.9}} \put(7.2,1.9){\vector(0,-1){0.9}}
\put(7.5,1.2){$\b$}
\put(5.1,0.8){\vector(1,0){1.9}} \put(6.9,0.8){\vector(-1,0){1.9}}
\put(5.9,0.4){$\a$}

\thicklines
\put(1,3){\line(1,1){4}}
\put(5,6){\line(1,1){1}}
\put(6,1){\line(1,1){1}}
\end{picture}
\caption{The map $\Tab$.}
\label{fig1}
\end{figure}
By identifying the points 0 and 1 we get an interval translation
map of the circle with two intervals. 
In the three-parameter class of all ITMs on the circle with two intervals,
the condition that $\lim_{y \to \beta_1^-} T(y) = \beta_3$ lets us
consider this two parameter subfamily.
A special example of an ITM of this form was considered in \cite{BK}.

\begin{proposition}\label{prop1}
If $\Tab$ is aperiodic then its restriction to $\Ob$ is  minimal.
If a single orbit of $\Tab$ is finite, then every orbit is eventually
periodic and the restriction to the attractor is isomorphic to a 
rational circle rotation.
\end{proposition}

\begin{proposition}\label{prop0}
The map $\Tab$ is of finite type if and only if  there is an interval $J$ such that
$\Tab$ induced on $J$ is an interval exchange of two intervals,
i.e., $\Tab$ is isomorphic to a circle rotation. If $\Tab$ is aperiodic and of
finite type, then it is uniquely ergodic.
\end{proposition}

We remark that if $\alpha =1$ or $\alpha = \beta$ then $T$ is a circle rotation,
while if $\beta = 0$ then $T$ is a noninvertible ITM on two intervals and 
in this case the identity on its attractor.

\begin{proofof}{Proposition \ref{prop0}}
We begin by proving the if statement.
By Proposition \ref{prop1}, if $T$ has a periodic point then it is of finite 
type.
Thus we assume that $T$ is aperiodic and therefore minimal on $\Ob$ by Proposition \ref{prop1}.
Let $J$ be an interval such that the induced map $T_J$
on $J$ is an interval exchange of two intervals. 
We have $J \subset \O$
and thus $T$ is of finite type since for any transitive ITM of infinite type
$\Omega$ must be a Cantor set \cite{ST}.

We prove the only if statement by contradiction.
Consider the first return map $\hT = T_{\Delta_1}$ of $\Tab$ to the interval
$\Delta_1 := [1-\a,1)$.  Clearly we have $\hT|_{[1-\a,1-\b)} = T|_{[1-\a,1-\b)}$.
There are two cases for the interval $[1-\b,1)$.  

i) The whole interval $[1-\b,1)$ returns at the same 
time.  In this case the map $\hT$ is an interval translation map on
two intervals.  Thus $\hT$ is an interval exchange of two intervals.

ii)  There are two different return times, that is, there is a positive
integer $k$ such that the left part of
the interval returns in $k$ steps and the right part of the interval
returns in $k-1$ steps.  In this case
the map $\hT$ of the form $T_{\a',\b'}$ where
\begin{equation}\label{Tab}
(\a',\b') \eqdef G(\a,\b) \eqdef \Big ( \frac{\b}{\a},
\ \frac{\b-1}{\a} + \Big \lfloor \frac{1}{\a} \Big \rfloor \Big ),
\end{equation}
where $k = \lfloor \frac{1}{\alpha} \rfloor$.
 
Repeating the inducing procedure we either are always in case ii) or
at some time step we reach case i), and then there is a subinterval
$J$ such that the first return map on $J$ is an interval exchange of two intervals.
If we are always in case ii) then the induced map to $\Delta_n := [1 - \prod_{i=0}^{n-1} \alpha_i,1)$
(where $(\a_i,\b_i) \eqdef G^i(\a,\b)$) is not invertible. This means that for each $n$ the set
$\Delta_n$ is not contained in $\O$.  
Since $\O \cap \Delta_n = \cup_{i \ge 0} T_{\Delta_n}^i(
\Omega \cap \Delta_{n+1})$ this implies that $\O$ cannot contain any interval
and thus $T$ is of infinite type.

Aperiodic circle rotations are always uniquely ergodic.  Since the Kakutani
tower of $\Tab$ over the induced map is finite, the map $\Tab$ is
uniquely ergodic if and only if the induced rotation is uniquely ergodic.
\end{proofof}

\begin{figure}[ht]
\setlength{\unitlength}{1cm}
 \begin{picture}(10,6)
 \put(0,0){\resizebox{2.7 cm}{5 cm}{\includegraphics{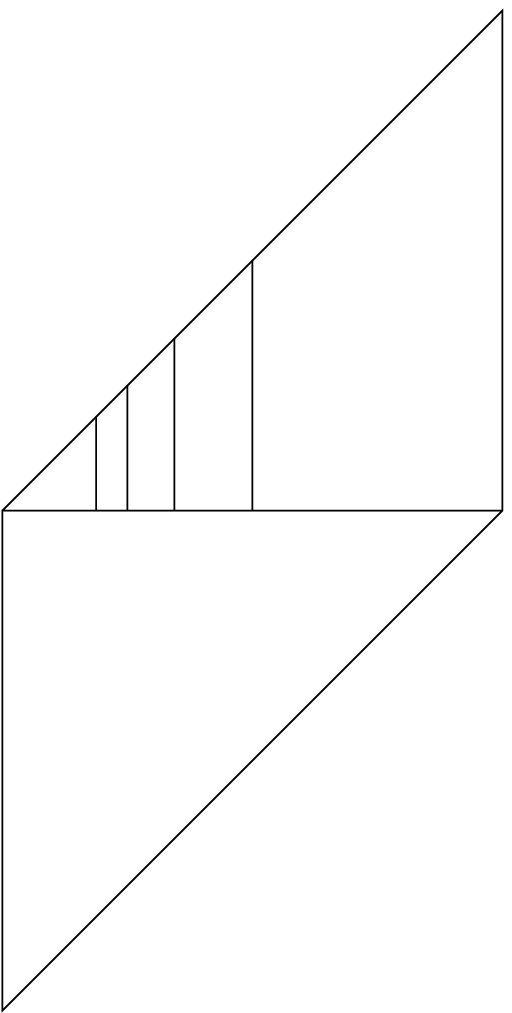}}}
 \put(5,0){\resizebox{6 cm}{5 cm}{\includegraphics{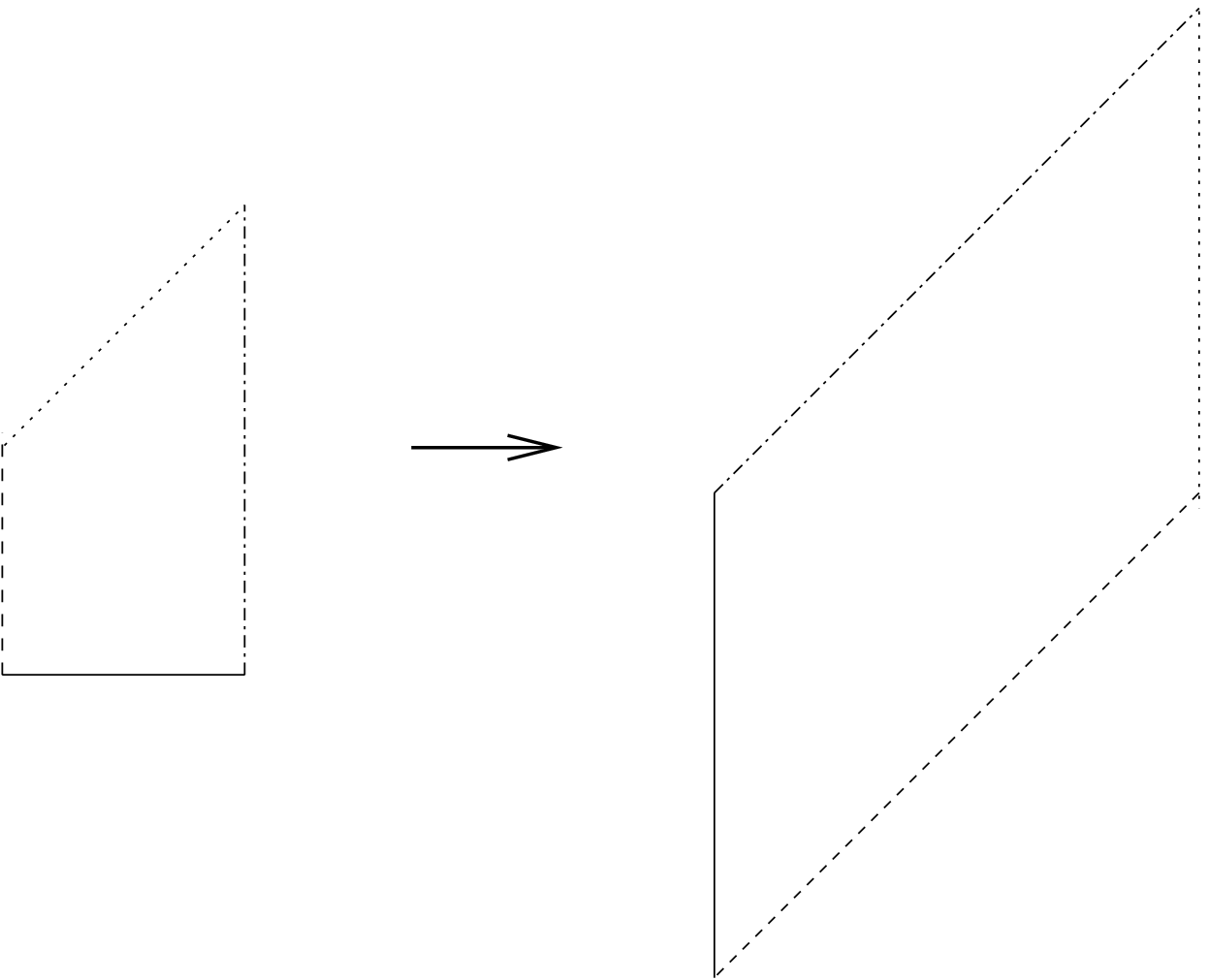}}}
\put(1.8,3.2){$\U_1$}
\put(.8,2.95){${\phantom{1}}_{\U_2}$}
\put(1,1.6){$L$} 
\put(0.3,2.65){.}
\put(0.25,2.6){.}
\put(0.2,2.55){.}
\put(7.15,3){$G$}
\put(5.5,2.3){$\U_k$}
\put(8.85,2.35){$R = U \cup L$}
\end{picture}
\caption{The action of the map $G$.}
\label{figextra}
\end{figure}

\begin{figure}[ht]
\setlength{\unitlength}{1cm}
 \begin{picture}(10,8)
 \put(1,-4){\resizebox{8 cm}{16 cm}{\includegraphics{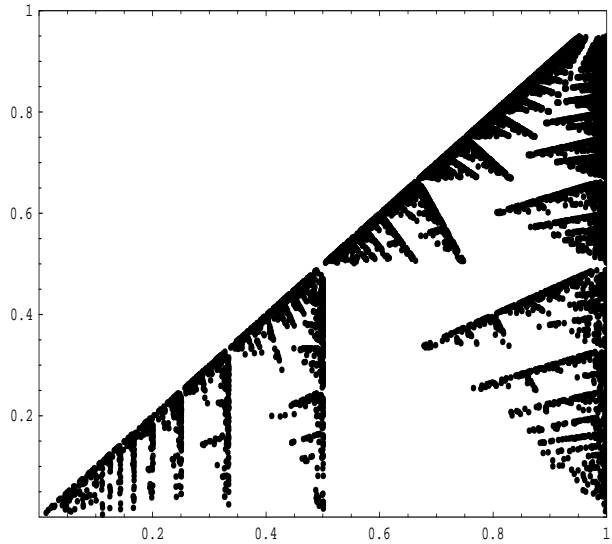}}}
 \end{picture}
\caption{Approximation of the set $A$ based on an iterated function system
with $10,000$ pixels.  The black should reach up to the top right corner
$(1,1)$, but does not because this is a neutral fixed point of $G$ 
(see Proposition \ref{propExpan}). 
}
\label{fig2}
\end{figure}

Equation \eqref{Tab} defines the {\em Gauss map} $G: U^{\circ} \to R$. To understand the action of this
map consider the region
$\U_k \eqdef \{ \ab \in U: \frac1{k+1} < \a \leq \frac1k\}$, see Figure~\ref{figextra}.
For each $k \in \N$, $G$ maps $\U_k$ into $R$ in the following
way:  the right boundary of $\U_k$ is mapped onto the top boundary of $R$, 
the left boundary of $\U_k$ is mapped onto the bottom boundary of $R$,
the bottom boundary of $\U_k$ is mapped onto the left boundary of $R$ 
and the top boundary of $\U_k$ is mapped onto the right boundary of $R$.
Note that if $\ab \in \partial \U_k$ for some $k$, then $T_{G(\a,\b)} \in \partial(R)$,
thus of finite type and consequently $\Tab$ is also of 
finite type. Let 
$$A \eqdef \cap_{n \ge 0} G^{-n}(U^{\circ}),$$  see Figure~\ref{fig2}.

\begin{corollary}
$\Tab$ is of infinite type if and only if $\ab \in A$.
\end{corollary}

\begin{proof}
If $(\a,\b) \not \in A$ then at some step of the inducing procedure
$G^n(\a,\b) \in L \cup \partial U$ and
the induced map is an ITM on two intervals.

If $(\a,\b) \in A$ then the inducing procedure can be repeated indefinitely
and thus $T$ is of infinite type by the proof of Proposition \ref{prop0}.
\end{proof}

\begin{theorem}\label{thmMeasA}
The set $A$ has Lebesgue measure $0$.
In particular, 
for Lebesgue almost every $\ab \in U$ the map $\Tab$ is of 
finite type and aperiodic.
\end{theorem}

Let $\NN = \{ 1,2,3,\dots \}$. The starting result of the analysis of 
our family of ITMs is:

\begin{theorem} \label{thmSub} The
set $A$ is uncountable and is naturally indexed 
by 
\[
\K \eqdef \{ k = k_0k_1\dots ; k_i \in \N, k_{2i} \neq 1 \mbox{ inf. often 
and } k_{2i +1} \neq 1 \mbox{ inf. often}
\}.
\] 
for any sequence $k$, the map $\Tab|_{\Ob}$ is isomorphic to the shift space
$(\Sab,\sigma)$ generated by the chain of substitutions 
$\chi_{k_0} \circ
\chi_{k_1} \circ \chi_{k_2} \circ \cdots$
where 
$$\chi_k \eqdef \left\{
\begin{array}{lll} 
1 & \to & 2\cr
2 & \to & 31^k\cr
3 & \to & 31^{k-1}
\end{array} \right.
$$
\end{theorem}

We will call an ITM 
{\em self-similar} if there exists one and hence infinitely
many subintervals
$\Delta \neq I$ such that the induced mapping $T_{\Delta} : \Delta
\to \Delta$ differs from $T$ only by affine scaling.

\begin{proposition}\label{propPeriodic}
The set of periodic points of $G$ is countable and dense in $A$. 
The set $\Omega$ is self-similar for any $\ab$ which
is $G$--periodic.
\end{proposition}

For each $x \in \O$ define the {\em itinerary} 
$e(x) \eqdef e_0(x) e_1(x) \dots$ where 
$e_k(x) = j$ iff $T^k(x) \in B_j = [\beta_{j-1},\beta_j)$.  
Let $\Sigma = \overline{e(\Omega)}$. 
The subshift $\Sigma$ is called {\em adic} if it is generated by
a sequence $l$, where 
$l = \lim_n \sigma_1 \circ \sigma_2 \circ \dots \circ \sigma_n(a)$,
and the substitutions $\sigma_i$ are primitive and
taken from a finite collection. 

\begin{proposition}\label{propAdic}
Let $\ab \in A$.  If the $G$-orbit of $\ab$ does not 
accumulate at the points $(0,0)$ or $(1,1)$ then
$\Tab$ is adic. 
This set of $\ab$ is uncountable and dense in $A$.
\end{proposition}
Write $\hd$, $\ubd$ and $\lbd$ for Hausdorff dimension and upper resp. 
lower box dimension. 
\begin{theorem}\label{thmHD}
There exists $r = 0.84955\cdots < 1$ such that for all $\ab \in A$, 
\[
\ubd(\Ob_{\a,\b}) \leq r.
\]
There exist maps $\Tab$ such that 
$\hd(\Ob_{\a,\b}) = \lbd(\Ob_{\a,\b}) = 0$,
as well as maps where 
$0 = \hd(\Ob_{\a,\b}) < \ubd(\Ob_{\a,\b})$.
\end{theorem}

Results of Boshernitzan imply that $\hd(\mu) \ge \frac12$ for any 
$T$-invariant Borel probability measure 
for almost every $\a,\b$ and for any $\a,\b$ which are algebraic \cite{Bos}.
The Hausdorff dimension of $\Ob$ has been computed before in
\cite{BK} for one particular example.
Let
\begin{equation}\label{Pk}
P_k(x) = x^3 - x^2 - k x + 1.
\end{equation} 
Denote the roots of $P_k$ by $\rl_k < r_k \leq \ru_k$.
One can check that $\rl_k \approx -\sqrt{k} + \frac12$,
$r_k \approx \frac1k$ and $\ru_k \approx \sqrt{k} + \frac12$
for large $k$. 
If $\a_k = r_k$, then $(\a_k,\a_k^2)$ is a fixed point of $G$,
and the Hausdorff dimension of 
$\Ob_{ \a_k,\a_k^2 }$ equals $-\frac{ \log(\ru_k) }{\log r_k}$.
For the case $k = 3$, the details have been worked out in \cite{BK}.
If $k \to \infty$, then 
$\hd( \Ob_{\a_k,\a_k^2} ) \to \frac12$.

As we have seen in Proposition~\ref{prop0}, every $\Tab$ of finite type
is isomorphic to a circle rotation, and therefore uniquely ergodic.
If $\Tab$ satisfies the hypotheses of Proposition~\ref{propAdic} 
then unique ergodicity is relatively easy to prove.
Fluctuations combined with large values in the sequence $(k_i)$ 
complicate a general analysis. Both cases occur:

\begin{theorem}\label{thmNUE}
Let $\Tab$ be of type
infinity, with sequence $(k_i)$ as in Theorem~\ref{thmSub}.
If for some $\lambda > 1$,  
$k_i \geq \lambda k_{i-1}$ for all $i$ sufficiently large,
then $\Tab$ is not uniquely ergodic.
\end{theorem}

\begin{theorem}\label{thmUE}
Let $\Tab$ be of type
infinity, with sequence $(k_i)$ as in Theorem~\ref{thmSub}.  
Let $L_{2i} = \min\{ r \geq 1 ; k_{2i+r} \neq 1 \}$. If 
\begin{equation}\label{(a)} 
\sum_i \frac{k_{2i}-1}{k_{2i}} \sqrt{ \frac{1}{k_{2i-1} L_{2i}} } = \infty,
\end{equation}
or
\begin{equation}\label{(b)} 
\prod_{i \geq 1} 
\frac{k_{2i}}{k_{2i-1} + \frac1{L_{2i}}} = 0,
\end{equation}
(or either condition holds with $2i$ replaced by $2i-1$),
then $\Tab$ is uniquely ergodic.
\end{theorem}
Conditions \eqref{(a)} and \eqref{(b)} 
have a non-empty symmetric difference. For instance, if
$k_{2i} = 2^i$, $k_{2i-1} = 3^i$, then \eqref{(b)} applies and 
not \eqref{(a)}, while
the case $k_i = i+1$ 
is covered only by Condition \eqref{(a)}.
Theorem~\ref{thmUE} allows a corollary on the abundance of uniquely ergodic
ITMs in our class.
\begin{corollary}\label{coroGdUE}
The set $A_{UE} := \{ \ab \in A \ | \ \Tab \mbox{ is uniquely ergodic} \}$
is a dense $G_{\delta}$ set in $A$, i.e., for each compact set 
$K \subset A$, $A_{UE} \cap K$ is dense $G_{\delta}$. 
\end{corollary}
The method for proving that the map $\Tab$ is uniquely ergodic or not
has the same flavor as the method that Keane \cite{Keane} used
for certain IET's. In fact, it gives also a way to estimate the number of
ergodic measures, and so we retrieve a result proven in more generality
by Buzzi and Hubert \cite{BH}.

\begin{corollary}\label{corNumbMeas}
Each map $\Tab$ has at most two ergodic invariant measures.
\end{corollary}

\section{Transitive implies minimal}

In this section we give a direct proof of Theorem~\ref{thmTrans}.

\begin{proofof}{Theorem \ref{thmTrans}}
If $T$ is of finite type, i.e.,~if $T|_{\Omega}$ is an IET, minimality 
is well known consequence of transitivity (see for example
Corollary 14.5.11 in \cite{KH}).
So assume that $T$ is of type $\infty$.
Let $y \in \overline{\Omega}$ be such that $\overline{\Omega}$ is $\omega(y)$,
the omega limit set of $y$.

\medskip

{\bf Claim 1:} There exists $0 < i < r$ and an interval $J \subset I$ 
such that $J \cap \overline{\Omega} \ne \emptyset$ and
$J \cap \overline{\Omega} \subset \omega(\b_i)$.

Note that $\beta_0$ and $\beta_r$ are not discontinuity points of the map $T$
and therefore not necessary in Claim 1.

Start with $i= 1$ and $J_1 = I$. Clearly 
$J_1 \cap \overline{\Omega} \ne \emptyset$. 
If $J_1 \cap \overline{\Omega} \subset \omega(\b_1)$, 
then we are finished. 
Otherwise, there exists an interval $J_2 \subset J_1$,
intersecting $ \overline{\Omega}$, such that 
$\omega(\b_1) \cap J_2 = \emptyset$.
Next check if $J_2 \cap  \overline{\Omega} \subset \omega(\b_2)$.
If so, then we are finished. 
Otherwise, there exists an interval $J_3 \subset J_2$,
intersecting $ \overline{\Omega}$, such that 
$\omega(\b_2) \cap J_3 = \emptyset$.
Continue this way. Since there are only finitely many discontinuity points,
we arrive at some interval $J_s$, intersecting $ \overline{\Omega}$ but disjoint from
$\omega(\b_i)$ for all $0 \leq i < r$. But this implies that 
$J_s \subset \Omega_n$ for all $n$, and $T$ cannot be of type $\infty$.

\medskip

{\bf Claim 2:} There exists $0 < i < r$ such that 
$\overline{\Omega} = \omega(\b_i)$.

Let $i$ and $J$ be from Claim 1. Since $\orb(y)$ is dense in
$\overline{\Omega}$, there exists $T^k(y) \in J$. Therefore
$\omega(\b_i) \supset \omega(T^k(y)) \supset \overline{\Omega}$.

{\bf Claim 3:} $\overline{\Omega} = \omega(\b_i)$ for  each $0 < i < r$.
  
Let $i_0$ be the $i$ from Claim 2.
Let $X_i = \{ \b_j; 0 < j < r \mbox{ and } \b_j \in \omega(\b_i) \}$.
Clearly, $X_{i_0} = \{ \b_1, \dots, \b_{r-1}\} \cap  \overline{\Omega}$.
Note that if $\b_j \in X_i$, then $X_j \subset X_i$, and also
that $X_i \neq \emptyset$. Indeed, if $X_i = \emptyset$, then there exists
a neighborhood $U$ of $T(\b_i)$ such that 
$T^k(U) \not\owns \b_j$ for any $0 < j < r$ and any $k \geq 0$.
This implies that any $\Omega_n$ contains an interval of length $\geq |U|$,
contradicting that $T$ is of type $\infty$.

To prove the claim, we need to show $X_i = X_{i_0}$ for each $i$.
Assume by contradiction that $X_{i_1}$ is the/a largest set strictly
smaller than $X_{i_0}$.
Find $\eps > 0$ such that $|T^k(\b_i) - \b_j| < \eps$
implies $\b_j \in X_i$.
Find $\b \in \orb(\b_{i_1})$ and $z \in \orb(y)$ such that
$|\b - z| < \eps/2$ and iterate these points.
Whenever some $\b_j \in (T^k(\b), T^k(z))$ 
(in particular, $\b_j \in X_{i_1}$), continue iterating 
with $T^k(z)$ and $\b_j$. We get that each point
in $\orb(z)$ is no more than $\eps/2$ away from some iterate
of some point $\b_j \in X_{i_1}$.
Since $\orb(z)$ approximates every point in $X_{i_0}$ arbitrarily
closely, it follows that $X_{i_1} = X_{i_0}$.

Now we can finish the proof. Take $x \in \Ob$ arbitrary.
No neighborhood $U \owns T(x)$ can be iterated indefinitely
without being cut, or otherwise $T$ is not of type $\infty$.
Therefore there exists $0 < j < r$ such that $\b_j \in \omega(x)$.
By Claim 3, $\omega(x) \supset \overline{\Omega}$.
\end{proofof}

\begin{proofof}{Proposition \ref{prop1}}
Theorem 2.4 of \cite{ST} states that $\Tab|_{\Ob}$ is minimal when it is aperiodic.

If there is a finite orbit, then there is an interval of periodic orbits.
Let $n$ be the period.  We consider $T$ as a mapping of the circle, thus
the point $1-\b$ is a point of continuity of $T$. Let
$[a,b)$ a maximal interval such that $T^i$ is continuous on $[a,b)$ for 
$i = 1,2,\dots n$ and $T^n$ is the identity.
Let $B_0 := \cup_{i=0}^{n-1} T^i [a,b)$, the set $B_0$ consists of 
$n$--periodic points.  By maximality, the points $1-\a$ and $1$ must be 
right end points of two of the intervals which make up $B_0$.  
Thus, it follows that $T$ has no other periodic points.

Let $B := \cup_{i=0}^\infty T^{-i} B_0$ and $ C:= I \setminus B$.
The set $B$ is the basin of attraction of the periodic component $A$;
it consists of a countable union of half open intervals. The set
$C$ is $T$--invariant.

Suppose that the second statement in the proposition is not true,
i.e., the set $C$ is non-empty.  

The proof of this proposition is now a modification of the proof
of Theorem 2.4 of \cite{ST}. This theorem gives a sharp upper bound
on the number of minimal sets for an aperiodic ITM.
We remind the reader of the notation.  Let $D := \{\b_0,\dots,\b_{r-1}\}$
and for an interval $J$ let $D(J) := \{ \b_i \in D: \ T^n \b_i \in 
\inte{(J)} \text{ for some } n \ge 0\}$. We called an interval $J$ 
{\em a minimal interval} if
$D(J) \ne \emptyset$ and $D(J) = D(J_1)$ for all subintervals $J_1$ of $J$.
In the proof of Theorem 2.4 it was shown that under the assumption of
aperiodicity minimal intervals induce a partition of $D$ and that
there is no minimal interval $J$ for which $D(J)$ consists of a single $\b_i$.

If we drop
the assumption of aperiodicity, then the same proof shows that
minimal intervals induce a partition of those $\b_i \in D$ whose orbit
is not periodic or preperiodic. Furthermore, in our case since $1-\b$ is
a point of continuity we can assume $D := \{0,1-\a\}$. 
However, one of these two points must be the left endpoint of one of the 
interval making up $B_0$, thus it has a periodic orbit. 
Thus no minimal interval can exist
since for every minimal interval $D(J)$ consists of at least two points.
\end{proofof}

\section{Properties of the Gauss map (Rauzy induction)}

In this section we collect some metric properties of the map $G$
and describe the structure of the set $A$.

\begin{proposition}\label{propExpan}
The second iterate of the map $G$ is uniformly expanding outside 
any neighborhood of the neutral fixed point $(1,1)$.
\end{proposition}

\begin{proof}
The map $G$ is infinite-to-one.
Away from its singularities the Jacobian of $G$ is
\[
DG = \left (
\begin{matrix} 
\frac{-\b}{\a^2} & \frac1{\a}\\
\frac{1-\b}{\a^2} & \frac1{\a}
\end{matrix}
\right )
\]
with determinant $-\a^{-3}$. The eigenvalues of $DG$ are
\[
\lambda^{\pm} \eqdef \frac1{\a} \cdot \left (\frac{\a -\b}{2\a} \pm 
\sqrt{\left( \frac{\a -\b}{2\a}\right )^2 + \frac1{\a}}\right ).
\]
The point $(1,1)$ is the only nonexpanding 
fixed point of $G$: it is elliptic with
eigenvalues $\pm 1$.  The eigenvalue $+1$ does not occur at any other
point in $R$ (it occurs on the curve $\b = \a^{-1}+\a-\a^2$), 
while the eigenvalue $-1$ occurs exactly on the intersection
of the curve $\b = \a^2 + \a -\a^{-1}$ with $R$.  Because
of the influence of this curve we consider the map $G^2(\ab) = (\a'',\b'')$
where
$$(\a'',\b'') = 
\Big ( \frac{\a}{\b}\cdot \big (1+ \Big \lfloor \frac{1-\a}{\a} \Big \rfloor
\big ) + \frac{\b-1}{\b},
\frac{2\b-1}{\b} + \frac{\a}{\b}\Big \lfloor \frac{1-\a}{\a} \Big \rfloor
+
\Big \lfloor \frac{\a-\b}{\b} \Big \rfloor \Big ).$$
The Jacobian of $G^2$ is
$$DG^2 = \left (
\begin{matrix} 
\frac1{\b}\big (1+ \Big \lfloor \frac{1-\a}{\a} \Big \rfloor\big ) &
\frac{-\a}{\b^2}\big (1+ \Big \lfloor \frac{1-\a}{\a} \Big \rfloor\big )+
\frac1{\b^2}\\
\frac1{\b} \Big \lfloor \frac{1-\a}{\a} \Big \rfloor &
\frac1{\b^2} + \frac{-\a}{\b^2} \Big \lfloor \frac{1-\a}{\a} \Big \rfloor
\end{matrix}
\right )
$$
and the determinant is $\b^{-3}$.
Recall that $k = 1+ \lfloor \frac{1-\a}{\a} \rfloor$.
The eigenvalue of $DG^2$ are
\[
\Lambda^\pm \eqdef \frac{1+\a(1-k) + \b k \pm
\sqrt{(1+\a(1-k)+\b k)^2 - 4\b}}{2\b^2}.
\]
For $k = 1$, we get $\Lambda^+ = 1/\b^2$ and $\Lambda^- = 1/\b$.
Note also that $\Lambda^- = 1/\a$ on the line $\a = \b$.
For the general case,
the equation $\Lambda^- > 1$ is equivalent to
\[
\left[ (1+\a(1-k) + \b k) - 2\b^2 \right]^2 > (1-\a(1-k)+\b k)^2 - 4\b,
\]
which follows from $1 + \b^3 > \b(1-\a(1-k)+\b k)$.
Since $\b \leq \a \leq 1/k$ this is easily checked to be true
for all $\b \in (0,1)$.
Hence $G^2$ is hyperbolic expanding outside a neighborhood
of the fixed point $(1,1)$.
\end{proof}

\begin{proofof}{Theorem \ref{thmMeasA}}
Let $V \eqdef U^{\circ} \cap G^{-1} U^{\circ}$.
Clearly $A = \cap_{n \ge 0} G^{-2n}(V)$.
Fix a small neighborhood $N_{\e}$ of the point $(1,1)$. Consider
the first return map $\Fe$ to the set $U \setminus \Ne.$  Consider 
$\Ae \eqdef \cap_{n \ge 0} \Fe^{-2n} (U_1 \setminus \Ne)$.
We have $A = (1,1) \cup \bigcup_{\e > 0} A_{\e}$ since
the elliptic point $(1,1)$ is weakly repelling for $G$.
Since $\Fe^2$ is uniformly expanding the set $\Ae$ has zero measure,
and thus $A$ has zero measure as well.

It immediately follows that $\Tab$ is of finite type for
Lebesgue a.e.~$\ab$. However, among the finite type parameters, the
(eventually) periodic ones form a countable union of smooth curves,
so $\Tab$ is aperiodic for a.e.~$\ab$.
\end{proofof}

\begin{proofof}{Theorem~\ref{thmSub}}
We can give symbolic dynamics for the set $A$.
The sets $\U_k$ from Figure \ref{figextra} form a Markov partition for the map $G|_A$.
If $\ab \in \partial \U_k$ for some $k$, then $\Tab$ is of 
finite type. Indeed, if $(\a,\b) \in \partial U_k$ is on the bottom or left boundary (excepting
the points $(0,0)$ and $(1,0)$) then $G(\a,\b) \not \in U$ so the code of $(\a,\b)$ is
not defined.  On the other hand if $(\a,\b) \in \partial U$ is on the top or right
boundary then
$G(\a,\b)$ belongs to either the upper boundary or the right 
boundary of the triangle $U$ and $G$ permutes these two sets.
Hence the code of $\ab$ will eventually have a $1$ at every other position.
It follows that the image of $A$ under the coding map is indeed $\K$,
and hence uncountable.

If $\ab \in A \cap \U_j$ then $G \ab$ can be interpreted
symbolically as the substitution $\chi_j$ in the following sense:
Let $T$ be the ITM and $F$ be the first return map
to the interval $[1-\a,1)$.
\begin{itemize}
\item The left branch of $F$ is identical to the middle 
branch of $T$. We denote this symbolically as
$1 \to 2$. 
\item Assuming that $\frac1{k+1} < \a \leq \frac1k$
we get that the middle branch of $F$ involves one iterate
of the third branch of $T$ followed by $k$ applications
of the left branch of $T$: $2 \to 31^k$. 
\item Still assuming that $\frac1{k+1} < \a \leq \frac1k$
we get that the middle branch of $F$ involves one iterate
of the third branch of $T$ followed by $k-1$ applications
of the left branch of $T$: $3 \to 31^{k-1}$. 
\end{itemize}

Take $k_n = j$ if $G^n(\ab) \in \U_j$. Let $\Delta_n$ be the $n$-th inducing interval,
and assume we want to code orbits of the first return map $T_{\Delta_n}$
according to the natural partition into branch domains, using symbols
$\{ 1, \ 2,\ 3 \}$.
Following the above pattern, we see that if $x \in \Delta_n$
has itinerary $u$ for $T_{\Delta_n}$, then $\chi_{k_{n-1}}(u)$ is the
itinerary of $x$ for $T_{\Delta_{n-1}}$.
By induction, the limit
\[ 
s = s_0s_1s_2\dots = 
\lim_{i \to \infty} \chi_{k_0} \circ \dots \circ \chi_{k_i}(3)
\]
gives the itinerary of the point $1$. Recall that, since $1 \notin [0,1)$,
we defined $T(1) = \lim_{x \to 1^-}T(x)$. Thus
$\Tab^{i-1}(T(1)) \in B_j$ if and only if $s_i = j$.
Let $\Sab = \overline{ \{ \sigma^i(s) : i \geq 0\} }$
be the corresponding shift space.
Then the above coding extends to a map $h:\Ob \to \Sab$
which is continuous and one-to-one, except possibly at the countable
set $\cup_i \Tab^{-i} (\{ 1-\a, 1-\b\})$.
\end{proofof}

\begin{proofof}{Proposition~\ref{propPeriodic}}
The map $G|_A$ is coded by the full shift on a countable alphabet. Thus
there are countably periodic codes.
Since $G|_A$ is non-uniformly expanding, see Proposition~\ref{propExpan}, 
each periodic code corresponds to a
single point and the periodic points lie dense in $A$. 

We turn to the self-similarity assertion. 
Let $(\a_i,\b_i) = G^i\ab$. Let $\Delta_n := [1-\prod_{i=0}^{n-1}\a_i,1).$
The map $G$ rescales the induced map
to be defined on a interval of length 1, thus $T$ and
$T_{\Delta_n}$ differ only by scaling for all
$k \in \N$ if and only if $\ab$ is a $G$--periodic point with period $n$.
\end{proofof}

\begin{proofof}{Proposition~\ref{propAdic}}
Let $O_k$ be the set of points $\ab \in A$ such that 
$G^i(\ab) \notin  \cup_{j > k} \U_j$
for all $i \geq 0$. 
The set $O_k$ is a subshift of finite type (in fact, the full shift
on $k$ symbols), thus it is uncountable.  Clearly the set $\cup_k  O_k$
is dense in $A$.
By Theorem~\ref{thmSub},
for $\ab \in O_k$, $\Tab$ is isomorphic to the shift space $\Sab$.
The substitutions $\chi_k$ are primitive for all $k \geq 2$.
By assumption, $(1,1)$ is not an accumulation point of $\orb(\a,\b)$,
so the substitution $\chi_1$ appears only in blocks of bounded 
length, and each ``block'' $\chi_1^n \circ \chi_r$, $2 \leq r \leq k$,
is primitive. Hence $\Sab$ is adic. Note that it is uniquely
ergodic as well, but this holds more generally, see Theorem~\ref{thmUE}.
\end{proofof}

\section{Dimension results}

\begin{proofof}{Theorem~\ref{thmHD}}
We will compute upper bounds for the box dimension by finding
suitable covers of $\Ob$.
Let $(\a_i,\b_i) = G^i(\a,\b)$ for $i \geq 0$.
The map $G^i$ computes the parameters of the ITM that results from the
first return map $T_{\Delta_i}$ to the 
interval $\Delta_i := [1-\pi_i,1)$ where 
$\pi_i := \prod_{j=0}^{i-1} \a_j$, so $\pi_0 = 1$. 
The map $T_{\Delta_i}$ consists of three branches,
the left branch of $T_{\Delta_i}$ is defined on an interval of length
$\pi_{i,1} = \pi_i (1-\a_i)$, the middle branch is defined on an interval
of length $\pi_{i,2} = \pi_i (\a_i-\b_i)$, while the length of
the right branch is $\pi_i \b_i$.
Also put $\pi_{i,3} = \pi_i \a_i$.
Let $k_i = \lfloor \frac1{\a_i} \rfloor$.
Then the intersection of $\Ob$ and the domain of the left branch 
can be covered by $k_i-1$ intervals of length 
$\pi_i \b_i = \pi_{i+1,3}$ and one interval of length 
$\pi_i (\beta_i - \alpha_i\beta_{i+1}) = \pi_{i+1,2}$.
The domain of the middle branch 
is an interval of length $\pi_i (\a_i - \b_i) = \pi_{i+1,1}$
and the interval $\Delta_{i+1}$
is the union an interval of length $\pi_i (\a_i - \b_i) = \pi_{i+1,1}$
and an interval of length $\pi_i \b_i = \pi_{i+1,3}$.
This is illustrated in Figure~\ref{fig3}.

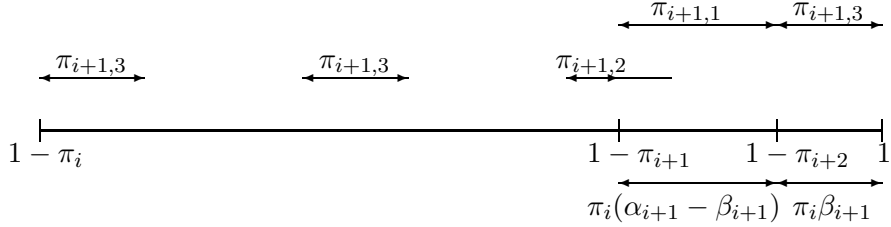
\begin{figure}[ht]
\unitlength=7mm
\begin{picture}(18,5)
\put(1,2){\line(1,0){16}}
\put(1,1.8){\line(0,1){0.4}}\put(0.4,1.4){$1-\pi_i$}
\put(12,1.8){\line(0,1){0.4}}\put(11.4,1.4){$1-\pi_{i+1}$}
\put(15,1.8){\line(0,1){0.4}}\put(14.4,1.4){$1-\pi_{i+2}$}
\put(17,1.8){\line(0,1){0.4}}\put(16.9,1.4){$1$}

\put(14,1){\vector(-1,0){2}} \put(13,1){\vector(1,0){2}}
\put(11.4,0.4){$\pi_i(\a_{i+1} - \b_{i+1})$}

\put(16.5,1){\vector(-1,0){1.5}} \put(15.5,1){\vector(1,0){1.5}}
\put(15.3,0.4){$\pi_i \b_{i+1}$}

\put(3,3){\vector(-1,0){2}} \put(1,3){\vector(1,0){2}}
\put(1.3,3.2){$\pi_{i+1,3}$}
\put(8,3){\vector(-1,0){2}} \put(6,3){\vector(1,0){2}}
\put(6.3,3.2){$\pi_{i+1,3}$}
\put(12,3){\vector(-1,0){1}} \put(11,3){\vector(1,0){1}}
\put(10.8,3.2){$\pi_{i+1,2}$}

\put(14,4){\vector(-1,0){2}} \put(13,4){\vector(1,0){2}}
\put(12.6,4.2){$\pi_{i+1,1}$}

\put(16.5, 4){\vector(-1,0){1.5}} \put(15.5, 4){\vector(1,0){1.5}}
\put(15.3, 4.2){$\pi_{i+1,3}$}

\thinlines
\put(12,3){\line(1,0){1}}
\end{picture}
\caption{Cover of $[1-\pi_i) \cap \Ob$ with intervals
$\pi_{i+1,j}$ for $k_{i+1} = 3$.}
\label{fig3}
\end{figure}

The map $T$ pushes the interval $\Delta_i$ and its subintervals
of length $\pi_{i,j}$ around through $[0,1)$.
The first return map to any image $T^k(\Delta_i)$ has the same structure
as $T_{\Delta_i}$, and hence $T^k(\Delta_i) \cap \Ob$ is covered by
the same number of intervals of length $\pi_{i,j}$ as $\Delta_i$.
Similarly, any interval of length $\pi_{i,j}$ can be replaced by a number
of intervals of length $\pi_{i+1,j'}$, according to the scheme:
\begin{eqnarray*}
\pi_{i,1} &\mbox{ is covered by }& 1 \times \pi_{i+1,2} \mbox{ and }
(k_{i+1}-1) \times \pi_{i+1,3}; \\ 
\pi_{i,2} &\mbox{ is covered by }& 1 \times \pi_{i+1,1}; \\
\pi_{i,3} &\mbox{ is covered by }& 1 \times \pi_{i+1,1} \mbox{ and }
1 \times \pi_{i+1,3}.  
\end{eqnarray*}
If we denote the number of such intervals by $l_{i,j}$, then the 
computation of their increase can be performed using a matrix:
\begin{equation}\label{cover}
\left( \begin{array}{c} l_{i+1,1} \\ l_{i+1,2} \\ l_{i+1,3} \end{array} \right)
=
\left( \begin{array}{ccc} 
0 & 1 & 1 \\ 1 & 0 & 0 \\ k_{i+1}-1 & 0 & 1 \end{array} \right)
\,
\left( \begin{array}{c} l_{i,1} \\ l_{i,2} \\ l_{i,3} \end{array} \right).
\end{equation}
The upper box dimension of $\Ob$ is bounded by
\[
\ubd(\Ob) \leq \limsup_i 
\frac{ \log(l_{i,1} + l_{i,2} + l_{i,3})}{-\log \pi_i}.
\]

Let $\rho = \log(\ru_2)/\log 2 = 0.84955\cdots$, where $\ru_2$ is the
leading root of $P_2$ defined by Equation \eqref{Pk}.
One can verify that $\ru_k \leq k^{\rho}$ for all $k \geq 1$.
Therefore we get for arbitrary $\ab \in A$:
\begin{eqnarray*}
\ubd(\Ob_{\a,\b})
&\leq&
\limsup_i \frac{ \log (l_{i,1} + l_{i,2} + l_{i,3}) }{ -\log \pi_i } \\
&\leq& \limsup_i 
\frac{ \log ( C \prod_{j=0}^i \ru_{k_j} ) }{ \log \prod_{j=0}^i k_j } \\
&=& \limsup_i \frac{ \sum_{j=0}^i \log \ru_{k_j} }
{ \sum_{j=0}^i \log k_j }  \leq \rho.
\end{eqnarray*}
This proves the upper bound.

As an example, if $\ab$ is the fixed point of $G$ in $\U_k$, $k\geq 2$,
then $\b = \a^2$ and $\a$ is the root $r_k \in (\frac1{k+1},\frac1k]$ 
of the polynomial $P_k$ from (\ref{Pk}).
The characteristic polynomial of the matrix in \eqref{cover} 
is also $P_k$ for $k = k_i$; recall that $\ru_k$ is its leading root.
Then we find $l_{i,j} = \OO(\ru_k^i)$ for each $j$, and 
therefore
$\ubd(\Ob_{\a,\b}) \leq -\log(\ru_k)/\log(r_k)$.
For $k = 3$ this bound 
equals $0.6635\cdots$ and in \cite{BK} (due to self-similarity) 
it is shown that if $k_i \equiv 3$, it is exact and equal to the Hausdorff dimension.
The same can be shown if $k_i \equiv k$ is constant.
If $k$ is taken very large, then $\ru_k \approx \sqrt{k} + \frac12$,
while $\a = r_k \approx \frac1k$.
Therefore 
$\frac12 \le \ubd(\Ob_{\a,\b}) \le \frac12 + \epsilon$ for 
sufficiently large $k$.
This is in agreement with the lower bound given in \cite{Bos}.
\medskip

Now to prove that $\hd(\Ob_{\a,\b}) = 0$ for some values of 
$\ab$ we argue as follows:
Let as before $\pi_i$ be the length of the domain of the $i$-th first return
map $T_{\Delta_i}$. 
We can cover $\Ob$ with $N_i := l_{i,1} + l_{i,2} + l_{i,3}$
intervals of length $\leq \pi_i$.
Given $N_i$ and $\pi_i$, choose $\a_i$ and $\b_i$
(thus determining $\b_{i+1}$) so that
\begin{equation}\label{hd0}
N_i \lfloor \frac1{\a_i} + 1 \rfloor (\pi_i \b_i)^{1/i}
+ N_i \lfloor \frac{\a_i}{\b_i} \rfloor 
(\pi_i \a_i {\b_{i+1}} )^{1/i} < 1. 
\end{equation}
This can be achieved as follows: write 
$\alpha_i = \frac{1}{k_i} - \e$ and 
$\beta_i=k_i\e + \e'$ for $0 < \e' \ll \e \ll (k_i)^{-2}$
and $0 < \e' < \alpha_i-k_i\e$.
Then (\ref{hd0}) follows from
\[
N_i \pi_i^{1/i} (k_i+1)^{1+1/i} \e^{1/i}
+ N_i \pi_{i+1}^{1/i} \frac{1}{k_i^2 \e} ( \e' )^{1/i} < 1, 
\]
which is easily satisfied for $(\e,\e')$ taken in an appropriate
region near $(0,0)$.

By the above reasoning, the domain of the left branch of $T_{\Delta_i}$ 
can be covered
by $\lfloor \frac1{\a_i} \rfloor$ intervals of length
$\pi_i \b_i$.
The domain of the middle branch of $T_{\Delta_i}$
is the domain of the left branch of $T_{\Delta_{i+1}}$ and can therefore
be covered by $\lfloor \frac1{\a_{i+1}} \rfloor
= \lfloor \frac{\a_i}{\b_i} \rfloor$ intervals of length
$\pi_i \a_i \b_{i+1}$.
Finally, the domain of the right branch of $T_{\Delta_i}$ is a single interval
of length $\pi_i \b_i$.
Putting these things together, we derive that $\Ob$ is covered
by $N_i \lfloor \frac1{\a_i} +1 \rfloor$ intervals of length
$\pi_i \b_i$ and $N_i \lfloor \frac{\a_i}{\b_i} \rfloor$
intervals of length $\pi_i \a_i \b_{i+1}$.
Due to the choice (\ref{hd0}),
the corresponding ``critical exponent'' of this cover is $\leq 1/i$.
If $\ab$ is such that indeed (\ref{hd0}) holds for infinitely many $i$, then
$\hd(\Ob_{\a,\b}) = \lbd(\Omega_{\a,\b}) = 0$.

In particular, if $(\a_i, \b_i)$ is a sequence which alternates 
values satisfying (\ref{hd0}) with long blocks of 
$k_i := \lfloor 1/\a_i \rfloor \equiv 3$, say,
then we will find that the upper box dimension 
$\ubd(\Ob_{\a,\b}) = 0.6635\cdots$
is larger than the Hausdorff dimension $\hd(\Ob_{\a,\b}) = 0$. 
This concludes the proof. 
\end{proofof}

\section{Invariant Measures}
In this section we prove Theorems~\ref{thmHMeas}, \ref{thmNUE}
and \ref{thmUE} as well as some related results.
Throughout the section, $\Tab$ is assumed to be of type infinity,
with corresponding sequence $(k_i)$.
Let us start with some notation.
Let $s = s_0s_1s_2s_3\dots$ 
be the fixed point of $\chi_{k_0} \circ \chi_{k_1} \circ \dots$.
For $a \in \{ 1,2,3 \}$, write
\[
\freq_a(\Tab) = \bigcap_{N \geq 1} \CH_{n \geq N, k \geq 1} 
\frac1n \# \{ k \leq i < k+n; s_i = a \},
\]
where $\CH$ denotes the convex hull over all $k \geq 1$ and $n \geq N$.
Without the subscript $a$, $\freq(\Tab)$ is the vector with three components.
We say that $\freq(\Tab)$ {\em exists} if each component
is a single point, which then satisfies
$\freq_1(\Tab) + \freq_2(\Tab) + \freq_3(\Tab) = 1$.

\begin{lemma}\label{freq}
If $\ab \neq (\hat\a,\hat\b)$, then
$\freq(\Tab) \cap \freq(T_{\hat\a, \hat\b}) = \emptyset$.
\end{lemma}

In other words, the frequency vector $\freq$ uniquely determines the
parameter $\ab$.

\begin{proof}
Each substitution $\chi_k$ has an associated matrix $A_k$
whose characteristic polynomial is $P_k$ as in (\ref{Pk}).
Define the simplex
$S = \{ (x,y,z) \in \RR^3 : 0\leq x,y,z, \, x+y+z = 1\}$;
the coordinates $x,y,z$ will play the role of
$\freq_1(\Tab)$, $\freq_2(\Tab)$, resp. $\freq_3(\Tab)$.
The matrix $A_k$ gives rise to a mapping $F_k$ on $S$ given by
\[ 
F_k(x,y,z) = \frac{1}{k(y+z)+x+y} (k(y+z)-z, x, y+z).
\]
Pass to new coordinates $\xi = x+y$ and $\eta = y+z$.
Since $x + y + z = 1$ on $S$, this gives a new map
\[
\tilde F_k(\xi,\eta) = (1- \frac{\eta}{k\eta + \xi}, \frac1{k \eta + \xi}),
\]
defined on the triangle 
$\tilde S = \{ (\xi,\eta) : \xi,\eta \leq 1, \, 1 \leq \xi+\eta\}$.
The map $\tilde F_k$ preserves lines, and the images
$\tilde F_k(\tilde S)$ are triangles ${\mathcal V}_k$ with corners
$(1,1)$, $(\frac{k-1}{k}, \frac{1}{k})$ and 
$(\frac{k}{k+1}, \frac{1}{k+1})$. These triangles
have disjoint interiors and tile the triangle $\tilde S$.
On ${\mathcal V}_k$, there is one inverse map
\begin{equation}\label{Finv}
\tilde F_k^{-1}(\xi,\eta) = \frac{1}{\eta}( 1+k(\xi-1), 1-\xi).
\end{equation}
The map $\tilde F_k$ has derivative
\[
D\tilde F_k(\xi,\eta) = 
\frac1{(k\eta + \xi)^2} \left( \begin{array}{cc}
\eta & -\xi \\ -1 & -k \end{array} \right),
\]
eigenvalues
$\lambda_{\pm} = \frac{1}{2(k \eta+ \xi)^2}
\left( k - \eta \pm \sqrt{ (k+\eta)^2+ 4\xi} \right)$,
and determinant $\frac{-1}{(k \eta+ \xi)^3}$. This determinant is less
than or equal to $1$ in absolute value. Write 
\[
Z_{i,j} = \tilde F_{k_j} \circ \tilde F_{k_{j+1}} \circ \dots
\circ \tilde F_{k_i}(\tilde S) \mbox{ and }
Z_j = \lim_{i \to \infty} Z_{i,j}.
\]
Since $\tilde F_k$ preserves lines but contracts area, $Z_{i,j}$ is
convex, while $Z_j$ is a point or a straight arc.

Since $(\a,\b) \neq (\hat\a,\hat\b)$ and $G^2$ is
(non-uniformly) expanding, see Proposition~\ref{propExpan},
$(k_i)$ and $(\hat k_i)$ are not the same.
Therefore $Z_{i,0}(\Tab)$ and $Z_{i,0}(T_{\hat\a,\hat\b})$ are not the 
same for some $i$. If they are disjoint, then 
also $Z_0(\Tab) \neq Z_0(T_{\hat\a,\hat\b})$.
Transforming back to the coordinates $x,y,z$, we obtain that
the frequency vectors $\freq(\Tab) \neq \freq(T_{\hat\a,\hat\b})$.

If  $Z_{i,0}(\Tab)$ and $Z_{i,0}(T_{\hat\a,\hat\b})$ meet in their
boundaries for all $i \geq 1$, then $k_j = \hat k_j = 1$
for all $i \geq i_0$. Note that $\tilde F_1$, under iteration,
contracts the triangle $\tilde S$ to the left upper corner $(0,1)$.
Therefore $Z_{i_0}(\Tab) =  Z_{i_0}(T_{\hat\a,\hat\b}) = (0,1)$.
Since $k_i \neq \hat k_i$ for some $i < i_0$, we get that
$Z_0(\Tab) \neq  Z_0(T_{\hat\a,\hat\b})$ after all.
\end{proof}

\begin{lemma}\label{lemmaUE}
Let $\Tab$ be of infinite type. Then
$Z_0(\Tab)$ is a single point if and only if
$\Tab$ is uniquely ergodic. 
\end{lemma}

It is clear that $\freq(\Tab)$ exists if and only if $Z_0(\Tab)$
is a single point. Moreover, $Z_0(\Tab)$ is a single point if and only if
$Z_i(\Tab)$ is a single point for some $i \geq 0$.

\begin{proof}
By Theorem~\ref{thmSub}, $(\Omega,\Tab)$
is isomorphic to a shift space $\Sab$ via an ``isomorphism''
which is one-to-one except on the countable set
$\cup_i \Tab^{-i}(\{ 1-\a, 1-\b\})$.
This set supports no invariant probability measure.
So it suffices to determine when $(\Sab,\sigma)$ is uniquely ergodic.

First assume that $Z_0$ is a singleton.
Let $\eps > 0$ be arbitrary and $C = c_0\dots c_M$ be any word.
Each $t \in \Sab$ is a concatenation of words 
of the form 
$W(i) = \chi_{k_0} \circ \dots \circ \chi_{k_N}(i)$, for $i \in \{ 1,2,3\}$
and some fixed $N$.
More precisely,
\[
t = WW(i_1)W(i_2)W(i_3) \dots
\] 
where $W$ is a suffix (possibly empty) of $W(i_0)$.
Let us say that an occurrence of $C$ in
$t$ {\em $N$-overlaps} if $t_{[k, k+M]} = C$ and
$k < |WW(i_1) \dots W(i_r)| < k+M$ for some $r$.
By taking $N$ sufficiently large, we can assume that
\[
\limsup_n \frac1n \#\{ i ; C = t_{[i, i+M]}
\mbox{ and $N$-overlaps}\} \leq \eps,
\]
uniformly over all $t \in \Sab$. Therefore
\[
\lim_n \frac1n 
\#\{ 0 \leq i < n :  C = t_{[i,i+M]} \mbox{ and does not $N$-overlap} \}
\]
differs from
$v(C,t) := \lim_n \frac1n \#\{ 0 \leq i < n ;  C = t_{[i,i+M]}\}$ 
by no more than $\eps$.

Since $Z_0(\Tab)$ is a singleton, there exists $\freq \in \tilde S$,
such that for each $\eps' > 0$ there is $M > N$ 
such that for each word
$V(b) = \chi_{k_{N+1}} \circ \dots \circ \chi_{k_M}(b)$,
\[
\left| \frac{\#\{ 0 < j \leq |V(b)| ; 
V(b)_j = a \} }{|V(b)|} - \freq_a \right| 
< \eps',
\]
for $a,b = 1,2,3$.
For each $a \in \{1,2,3\}$, let $C_a$ be the number of occurrences of $C$ 
in the word $W(a)$. 
Then
\[
v(C,t) = \frac{ \sum_{a=1}^3 C_a \cdot (\freq_a + \OO(\eps')) } 
{ \sum_{a=1}^3 |W(a)| (\freq_a + \OO(\eps')) } + \OO(\eps).
\]
Since $\eps$ and $\eps'$ are arbitrary, we see that
$v(C,t)$ is independent of the string $t$.
Thus unique ergodicity follows.

\medskip

Conversely, if $Z_0(\Tab)$ is not a single point,
it has diameter $\diam Z_0 := \delta > 0$.
It follows that $\diam Z_{i,0}(\Tab) \geq \delta$ for every $i \in \N$.
The extremal points of $Z_{i,0}$ are the images under
$\tilde F_{k_0} \circ \dots \circ \tilde F_{k_i}$ of the corners of $\tilde S$.
For simplicity, we can assume that the vertical
height of $Z_{i,0} \geq \delta/2$.   
Recall that the variables $x,y$ and $z$ give the frequencies of the 
symbols $1,2$ and $3$ in words of $\Sab$, and that $x = 1-\eta$.

There exists $a, b \in \{ 1,2,3\}$ 
such that for infinitely many $j$, we find that 
the frequencies of the symbol $1$:
\[
\frac{1}{|\chi_{k_0} \circ \chi_{k_1} \circ \dots \circ \chi_{k_j}(a)|}
\#\{ \mbox{symbols } 1 \mbox{ in }
\chi_{k_0} \circ \chi_{k_1} \circ \dots \chi_{k_i}(a) \}
\] 
and
\[
\frac{1}{|\chi_{k_0} \circ \chi_{k_1} \circ \dots \circ \chi_{k_i}(b)|}
\#\{ \mbox{symbols } b \mbox{ in }
\chi_{k_0} \circ \chi_{k_1} \circ \dots \chi_{k_i}(b) \}
\] 
differ by at least $\delta/2$. Hence $\Sab$ is not uniquely ergodic.
\end{proof}

\begin{proofof}{Theorem~\ref{thmNUE}}
According to Lemma~\ref{lemmaUE} we need to show that compositions
$\tilde F_{k_0} \circ \tilde F_{k_1} \circ \dots$ do not
contract the simplex $\tilde S$ to a single point.
The composition of two maps $\tilde F$ has the form
\[
\tilde F^2(\xi,\eta) = \tilde F_{k'} \circ F_k(\xi,\eta) 
= (1 - \frac{1}{k'+(k-1)\eta + \xi}, \frac{k\eta + \xi}{k'+(k-1)\eta) + \xi}).
\]
The map $\tilde F_{k_{2i-1}} \circ \tilde F_{k_{2i}}$ acts on the second
component as
\[
h_i : \eta \mapsto \frac{ k_{2i} \, \eta + \xi}{k_{2i-1} + 
(k_{2i}-1) \eta + \xi}.
\]
We will show
$\lim_{j\to\infty} h_1 \circ h_2 \circ \dots \circ h_j([0,1])$ is a 
non-degenerate interval.

Let $\tau = 1/(\lambda-1)$ and $\tau' = 1-1/\lambda$.
Take $i_0$ so large that $k_{i} \geq \lambda k_{i-1}$ for all
$i \geq i_0$ as well as $\tau < \tau' k_{2i_0-2}$.
Obviously, $h_i$ is an increasing M\"obius transformation on $[0,1]$.
Assume that $2i \geq i_0+2$. Then
\[
h_i(\tau') \geq \frac{\tau' k_{2i}}{k_{2i}/\lambda + \tau' (k_{2i}-1) + 1}
\geq \frac{\tau'}{1/\lambda + \tau'} = \tau'.
\]
It follows that $h_i([\tau',1]) \subset (\tau',1)$, and hence
\[
\lim_{j \to \infty} h_{i_0} \circ h_{i_0+1} \circ \dots \circ h_j(1) 
\geq \tau'.
\]
On the other hand, if $\eta \leq \tau/k_{2i}$, then 
$h_i(\eta) \leq \frac{1+\tau}{k_{2i-1}} \leq \frac{1+\tau}{\lambda k_{2i-2} } =
\frac{\tau}{k_{2i-2}}$.
By induction we find 
$\lim_{j \to \infty} h_{i_0} \circ h_{i_0+1} \circ \dots \circ h_j(0) 
\leq \frac{\tau}{k_{2i_0-2}}$.
Therefore
\begin{eqnarray*}
[p,q] &:=&
\lim_{j \to \infty} h_{1} \circ h_{2} \circ \dots \circ h_j([0,1]) \\
&=&  h_{1} \circ h_{2} \circ \dots \circ h_{i_0-1}(
\lim_{j \to \infty} h_{i_0} \circ h_{i_0+1} \circ \dots \circ h_j([0,1])) \\
&\supset&
h_{1} \circ h_{2} \circ \dots \circ h_{i_0-1}([\frac{\tau}{k_{2i_0-2}},\tau']).
\end{eqnarray*}
In particular, $p < q$.
Lemma~\ref{lemmaUE} implies that
$\Sab$ and therefore $\Tab$ are not uniquely ergodic.
\end{proofof}

\begin{proofof}{Theorem~\ref{thmUE}}
According to Lemma~\ref{lemmaUE}, we need to show that $Z_0$ 
(or equivalently $Z_1$) is a single point.
The first component of $\tilde F^2(\xi,\eta)$ contracts the
interval to a single point because the derivative 
with respect to $\xi$ is $\leq 1$ with equality
only if $k = k' = 1$ and $\xi = 0$. Hence $Z_0$ has 
``width'' $0$.
For the ``height'', we need the second component of 
$\tilde F^2(\xi,\eta)$, for which we will use the maps $h_i$ from 
the proof of Theorem~\ref{thmNUE}.
The compositions $H_i = h_1 \circ \dots \circ h_i$ are 
M\"obius transformations represented by the matrix-multiplications
\begin{eqnarray*}
B_i &=& \left( \begin{array}{cc}
b^i_{1,1} & b^i_{1,2} \\
b^i_{2,1} & b^i_{2,2} 
\end{array} \right) \\
&=& 
\left( \begin{array}{cc}
k_2 & \xi_1 \\
k_2-1 & k_1 + \xi_1 
\end{array} \right) 
\dots \left( \begin{array}{cc}
k_{2i} & \xi_i \\
k_{2i}-1 & k_{2i-1} + \xi_i 
\end{array} \right),
\end{eqnarray*}
where the $\xi_i \in [0,1]$. Since 
\[
H_i([0,1]) \subset [H_i(0), \lim_{\eta \to \infty} H_i(\eta)]
= [b^i_{1,2} / b^i_{2,2} \ , \ b^i_{1,1} / b^i_{2,1}],
\]
it suffices to show that $B_i$ contracts the cone
${\mathcal C} := (\RR_{\geq 0})^2$ to a one-dimensional subcone
as $i \to \infty$.
In order to do this, we use {\em Hilbert metric} on $\mathcal C$:
Given $v,w \in \mathcal C$, define
\[
\Theta(v,w) = \log \left(
\frac{ \inf \{ \mu ; \mu v - w \in {\mathcal C}\} }
{ \sup\{ \lambda ; w - \lambda v \in {\mathcal C} \} } \right).
\]
In fact, $\Theta$ is a semi-metric, because $\Theta(v,w) = 0$ if and
only if $v$ is a multiple of $w$.
Let $T: {\mathcal C} \to {\mathcal C}$ be a linear map.
It is shown in e.g.~\cite{Birk} that
$\Theta(T(v), T(w)) \leq \tanh(D/4) \Theta(v,w)$ for
$D = \sup_{v',w' \in T( {\mathcal C} ) }\Theta(v',w')$.
In particular, $T$ is a contraction if $T$ maps  
$\partial {\mathcal C} \setminus \{ 0 \}$
into the interior of $\mathcal C$.

In our setting, the transformations $T$ are represented by
matrices of the form 
$\left( \begin{array}{cc}
k & \xi \\
k-1 & k' + \xi 
\end{array} \right)$, 
and we can easily check that
$D$ is assumed by taking $v' = (1,0)^{\top}$ and
$w' = (0,1)^{\top}$, so
$D = \log \frac{k(k'+\xi)}{(k-1)\xi}$.
Hence the contraction factor is
\begin{eqnarray*}
\tanh(\frac{D}{4}) &=& 
\frac{ 
\root{4}\of{ \frac{k(k'+\xi)} {(k-1)\xi} } - 
\root{4}\of{ \frac{(k-1)\xi}{k(k'+\xi)} } }
{ 
\root{4}\of{ \frac{k(k'+\xi)}{(k-1)\xi} } + 
\root{4}\of{ \frac{(k-1)\xi}{k(k'+\xi)} } } \\
&=&
\frac{ \sqrt{k(k'+\xi)} - \sqrt{(k-1)\xi} }
{ \sqrt{ k(k'+\xi) } + \sqrt{(k-1)\xi} } \\
&=&
1-2 \frac{ \sqrt{ k(k-1)(k'+\xi) \xi }  - (k-1)\xi/2 }
{ kk'+\xi } \\
&\leq&
1- 4 \frac{k-1}{k} \sqrt{ \frac{\xi}{k'} }.
\end{eqnarray*}
The variable $\xi$ is the result of iterating
$\tilde h:\xi \mapsto 1 - \frac{1}{k'+\xi + (k-1)\eta}$, the first component
of $\tilde F^2$.
The image $\tilde h(\xi) \geq \frac12$ unless $k = k' = 1$.
If $k = k' = 1$, then $\xi = 0$ is an indifferent attracting fixed
point, and $\tilde h^n(1) = 1/(n+1)$.
Therefore, in the above calculation, we can estimate 
$\xi \geq 1/L_{2i}$.
Hence, the height of $Z_0$ is less than
$\prod_i \tanh(D_i/4) \leq \exp\left( -4 \sum_i \frac{k_{2i}-1}{k_{2i}}
\sqrt{ \frac{1}{L_{2i} k_{2i-1}} } \right)$.
The assumption on $(k_i)$ gives that $Z_0$ is indeed a single point.
This finishes the proof of the theorem with condition \ref{(a)}.

For condition (\ref{(b)}), as in the proof of the part with 
condition (\ref{(a)}), 
the width of the $Z_1$ is $0$.
The compositions of M\"obius transformations
$H_i = h_1 \circ h_2 \circ \dots \circ h_i$ satisfy,
\[
|H_i([0,1])| = \sqrt{ H'_i(0) \cdot H'_i(1) } \leq
H'_i(0) \leq \prod_{j=1}^i h'_j(0),
\]
where the last inequality follows because each $h_i$
is increasing with decreasing derivative.
We compute $h'_j(\eta) = \frac{ k_{2j} k_{2j-1} + \xi}
{ (k_{2j-1} + (k_{2j}-1)\eta + \xi)^2}$,
and therefore $h'_j(0) \leq k_{2j}/(k_{2j-1}+1/L_{2j})$.
It follows that 
\[
\diam Z_1 \leq \lim_j H_i([0,1]) \leq \lim_i 
\prod_{1 \leq j \leq i} 
\frac{k_{2j}}{k_{2j-1} + \frac1{L_{2i}}} = 0.
\]
Lemma~\ref{lemmaUE} yields unique ergodicity.
\end{proofof}

\begin{proofof}{Corollary~\ref{coroGdUE}}
Note that the coding map $(\a,\b) \mapsto k \in \K$
is continuous on $A$.
Hence it suffices to consider the space $\K$.
For any cylinder 
$C_{e_1\dots e_n} = \{ k \in \K \ | \ k_i = e_i, \ i = 1, \dots n\}$, let 
$U^m_{e_1\dots e_n} = \{ k \in \K \ | \ k_i = e_i, \ i = 1, 
\dots n, k_{n+i} = 2, \ i = 1, \dots m \}$.
Clearly ${\mathcal U}^m = \cup_n \cup_{e_1 \dots e_n} U^m_{e_1\dots e_n}$
is open and dense in $\K$. Moreover, for each $k \in  U^m_{e_1\dots e_n}$
we have
\[
\sum_i \frac{k_{2i}-1}{k_{2i}} \sqrt{ \frac{1}{k_{2i-1} L_{2i}} } 
\geq \frac{m-2}{4\sqrt{2}}. 
\]
Therefore $\cap_m {\mathcal U}^m$ is a dense $G_{\delta}$ set of sequence
$k$ satisfying Condition~\eqref{(a)} 
This proves the corollary.
\end{proofof}

\begin{proofof}{Corollary~\ref{corNumbMeas}}
The maximum of two ergodic measures corresponds to the
at most two extremal points of the sets $Z_0$ in the proof of 
Lemma~\ref{freq}, see~\cite{Keane}.
This also follows immediately from the result of Buzzi and Hubert, \cite{BH}.
\end{proofof}

The next result gives a candidate (modulo finiteness) 
of an invariant measure.

\begin{proofof}{Theorem \ref{thmHMeas}}
Let $\V$ be a cover of $\Omega$ by the intervals forming $\Omega_n$
satisfying $J \cap \Omega \ne \emptyset$.
The map $T$ induces a multivalued map $\psi$ of the $J$'s by 
$\psi(J) = J'$ if $T(J) \cap J' \neq \emptyset$. Thus $\psi$
has one value at $J$ if $T|_J$ is continuous, otherwise it has two values
provided $n$ is sufficiently large.
As a result, at most $r$ of the $J's$ can have more than one $\psi$-preimages,
where $r$ is the number of discontinuity points.

Take $\e > 0$ so small that any two $J,J' \in \V$ are at least $\e$ apart. 
Let $\U \eqdef \{U_i\}$ be an open cover of $\Omega$ with the diameters
of the $U_i$ all less than $\e$.  
Suppose $J$ is an interval such 
$\psi(J)$ has exactly one preimage.
If the subcover of $\U$ covering $J$
gives a good approximation of $H_d(J)$, i.e., $\sum_{U_i \cap{J \ne \emptyset}}
\diam(U_i)^d \approx H_d(J)$, then the translated subcover 
$\{ T(U_i) \}_{U_i \cap J \neq \emptyset}$ satisfies
$\sum_{U_i \cap{J \ne \emptyset}} \diam(T(U_i))^d \approx H_d(T(J))$.
Since there are only finitely many 
intervals $J$ such 
$\psi(J)$ has more than one preimage, the union of these intervals
is negligible as the $n \to \infty$ (and hence $\e \to 0$).
So when minimizing over all $\e$-covers $\U$, we can restrict ourselves 
to $T$-invariant $\e$-covers and 
find that Hausdorff measure is $T$-invariant.
\end{proofof}

\section{Acknowledgements}
We thank Pascal Hubert for useful conversations and the anonymous referee for
many valuable suggestions.


\begin{thebibliography}{99}

\bibitem[AKT]{AKT} R.~Adler, B.~Kitchens, C.~Tresser, 
{\em Dynamics of non-ergodic piecewise affine maps of the torus,}
Ergod. Th. Dyn. Sys. {\bf 21} (2001) 959-1000.

\bibitem[ACP]{ACP} P.~Ashwin, W.~Chambers, G.~Petrov,
{\em Lossless digital flow oscillations; approximation of invariant fractals,}
Intl. J. Bifurcation and Chaos, {\bf 7} (1997) 2603--2610.


\bibitem[B]{Birk} G.~Birkhoff,
{\em Extensions of Jentzsch's theorem,}
Trans. Amer. Math. Soc. {\bf 85} (1957) 219--227.

\bibitem[Bo]{Bo1} M.~Boshernitzan, 
{\em A unique ergodicity of minimal symbolic flows with
linear block growth rate,} Journal d'Anal.~Math.~{\bf 44} (1984/5) 77--96.

\bibitem[Bo1]{Bos} M.~Boshernitzan, 
{\em Quantatitive recurrence results,}
Invent. Math. {\bf 113} (1993) 617--631.

\bibitem[BK]{BK} M.~Boshernitzan, I.~Kornfeld, 
{\em Interval translation mappings,}
Ergod.\ Th.\ Dyn.\ Sys.\ {\bf 15} (1995) 821--831.

\bibitem[BH]{BH} J.~Buzzi, P.~Hubert,
{\em Piecewise monotone  maps without periodic points: Rigidity, measures
and complexity,} IML preprint 2001--27.

\bibitem[C1]{Cas94b} J.~Cassaigne, {\em Facteurs sp\'eciaux et complexit\'e,}
Bull.~Belg.~Math.~Soc.~{\bf 4} (1997) 67--88.

\bibitem[C2]{Cas2003} J.~Cassaigne,
{\em Computing the subword complexity of an s-adic sequence:
application to a family of interval translation maps,}
in preparation.

\bibitem[G]{G} A.~Goetz, 
{\em Dynamics of a piecewise rotation,} 
Discrete Contin. Dynam. Systems {\bf 4} (1998) 593--608.

\bibitem[HR]{HR} F.\ Hofbauer, P.\ Raith, 
{\em Topologically transitive subsets of piecewise
monotonic maps, which contain no periodic points,} 
Monatsh. Math. {\bf 107} (1989) 217--239.

\bibitem[KH]{KH} A.~Katok, B.~Hasselblatt,
{\em Introduction to the modern theory of dynamical systems,}
Cambridge Univ.~Press (1995).

\bibitem[K]{Keane} M.~Keane,  
{\em Non-ergodic interval exchange transformations,} 
Israel J. Math. {\bf 26} (1977) 188--196. 

\bibitem[L]{L} G.~Levitt, {\em La dynamique des pseudo-groupes de
rotations,} Invent.~Math.~{\bf 113} (1993) 633--670.

\bibitem[LV]{LV} J.~Lowenstein, F.~Vivaldi, {\em 
Embedding dynamics for round-off errors near a periodic orbit,} 
Chaos {\bf 10} (2000) 747--755.

\bibitem[M]{M} H.~Masur, {\em 
Interval exchange transformations and measured foliations,} 
Ann. of Math. {\bf 115} (1982) 169--200. 

\bibitem[R]{R} G.~Rauzy, {\em Echanges d'intervalles et transformations 
induites,} Acta Arith.~{\bf 34} (1979) 315--328.

\bibitem[ST]{ST} J.~Schmeling, S.~Troubetzkoy, {\em Interval translation
mappings,} in ``Dynamical systems from crystals to chaos'' J.-M.~Gambaudo
et al.~eds. World Scientific, Singapore (2000) 291--302.

\bibitem[V1]{V1} W.~Veech, {\em Gauss measures for transformations on the 
space of interval exchange maps,} Ann.~Math.~{\bf 115} (1982) 201--242.

\bibitem[V2]{V2} W.~Veech, {\em The metric theory of interval exchange 
transformations I,II,III,} Amer.~J.~Math.~{\bf 106} (1984) 1331--1421.

\end{thebibliography}
\end{document}